\documentclass[11pt,letterpaper,reqno]{amsart}
\usepackage{tikz}
\usetikzlibrary{positioning, shapes.geometric, arrows.meta, calc}
\usepackage{amssymb}
\usepackage{amsmath}
\usepackage{amsthm}
\usepackage{amsfonts}
\usepackage{bbm}
\usepackage{enumitem}
\usepackage{xcolor}
\usepackage{pgfplots}
\pgfplotsset{compat=1.18} 
\usepackage{booktabs}
\usepackage{graphicx}
\usepackage[T1]{fontenc}
\usepackage{doi}
\usepackage{float} 
\usepackage{comment} 

\let\originalleft\left
\let\originalright\right
\renewcommand{\left}{\mathopen{}\mathclose\bgroup\originalleft}
\renewcommand{\right}{\aftergroup\egroup\originalright}

\newcommand{\ZZ}{\mathbb{Z}}

\newcommand\doublecheck{\textcolor{blue}{\checkmark\kern-0.5em\checkmark}}
\newcommand{\newvtheorem}[2]{\newtheorem{#1}[theorem]{\llap{\textnormal{\doublecheck}\ }#2}}

\theoremstyle{plain}

\newvtheorem{vtheorem}{Theorem}

\newvtheorem{vlemma}{Lemma}

\newvtheorem{vclaim}{Claim}

\theoremstyle{definition}

\theoremstyle{remark}
\newtheorem*{remark}{Remark}
\newvtheorem{vremark}{Remark}

\begin{document}
\title[Matching integers to distinct multiples]{Optimal bounds for an Erd\H{o}s problem on matching integers to distinct multiples}

\author[W.~van Doorn]{Wouter van Doorn}
\address{Groningen, the Netherlands} 
\email{wonterman1@hotmail.com}

\author[Y.~Li]{Yanyang Li}
\address{School of Mathematics, Southeast University, Nanjing 211189, P. R. China}
\email{liyanyang1219@gmail.com}

\author[Q.~Tang]{Quanyu Tang}
\address{School of Mathematics and Statistics, Xi'an Jiaotong University, Xi'an 710049, P. R. China}
\email{tang\_quanyu@163.com}



\begin{abstract}
Let $f(m)$ be the largest integer such that for every set $A = \{a_1 < \cdots < a_m\}$ of $m$ positive integers and every open interval $I$ of length $2a_m$, there exist at least $f(m)$ disjoint pairs $(a, b)$ with $a \in A$ dividing $b \in I$. Solving a problem of Erd\H{o}s, we determine $f(m)$ exactly, and show
$$
f(m)=\min\bigl(m,\lceil 2\sqrt{m}\,\rceil\bigr)
$$
for all $m$. The proof was obtained through an AI-assisted workflow: the proof strategy was first proposed by ChatGPT, and the detailed argument was subsequently made fully rigorous and formally verified in Lean by Aristotle. The exposition and final proofs presented here are entirely human-written.
\end{abstract}

\maketitle

\section{Introduction}\label{sec:intro}

For each positive integer \(m\), let \(f(m)\) denote the largest integer \(r\) with the following property: whenever
\[
A=\{a_1<a_2<\cdots<a_m\}
\]
is a set of \(m\) positive integers, then for every real number \(x\) there exist distinct elements \(c_1,\dots,c_r \in A\) and distinct integers \(b_1,\dots,b_r\)
such that
\[
x < b_i < x + 2a_m \qquad \text{and} \qquad c_i\mid b_i \qquad (1\le i\le r).
\]

In their 1959 paper \cite{ErdosSuranyi59}, Erd\H{o}s and Sur\'anyi proved the first result on $f(m)$, obtaining the lower bound $f(m) \ge \sqrt{m}$. In the opposite direction, Erd\H{o}s and Selfridge \cite{Erdos78, Erdos86} showed that $f(m^2) \le 2m$ holds for all $m$, from which one can deduce the general upper bound $f(m) \le 2 \lceil \sqrt{m}\, \rceil$. Hence, known lower and upper bounds on $f(m)$ differ by a factor of two, and it is natural to wonder what the exact growth rate is.
In \cite{Erdos95}, Erd\H{o}s in particular asked whether the lower bound can be improved, and estimating $f(m)$ is now recorded as Erd\H{o}s Problem \#650 on Bloom's website \cite{Bloom650}. In this paper, we determine \(f(m)\) exactly, thereby completely resolving the problem. 

Taking a step back, one might wonder why $2a_m$ would be a natural choice for the length of the interval in the definition of $f(m)$. For a start, if we instead let the length of the interval be $ca_m$ for some $1 < c < 2$, then the corresponding definition of $f(m)$ would actually give $f(m) = 1$ for all positive integers $m$. Indeed, such an interval certainly contains an integer divisible by $a_1$, giving the lower bound $f(m) \ge 1$. For the corresponding upper bound one can choose a large integer $M$, choose $A = \{a_1, a_2, \ldots, a_m\}$ with $a_i = M+i$ for all $i$, and set $x = -a_1 + \prod_{i=1}^{m} a_i$. If $M$ is sufficiently large in terms of $c$ and $m$, then one can verify that, with $1 < c < 2$, the only integer in $(x, x+ca_m)$ divisible by some element of $A$ is $x + a_1$. For $c \ge 3$ on the other hand, Erd\H{o}s writes in \cite{Erdos78} and \cite{Erdos86} that ``all hell breaks loose''. For such values of $c$ it is unclear what to expect exactly, although in \cite{Erdos86} a quick argument is recorded showing that any interval of length $3a_m$ contains more than $\sqrt{6m}$ distinct multiples of the $a_i$. Finally, even though we defined $f(m)$ using an interval length of $2a_m$, the results in this paper will actually cover the entire intermediate range $2 \le c < 3$. Indeed, in analogy with \cite{Erdos78, Erdos86}, it will become clear that our proofs still go through with $2$ increased to any real smaller than $3$. 

Now, to better understand \(f(m)\), it is convenient to reformulate its definition in graph-theoretic terms. Writing $\mathbb{N}=\{1,2,3,\dots\}$ as the set of positive integers, let a set $A\subseteq \mathbb{N}$ (with $|A| = m$ and $\max A = a_m$) and a real number $x$ be given. With $B=(x,x+2a_m)\cap\mathbb{Z}$, we then define a bipartite graph \(G(A,x)\) with vertex classes \(A\) and \(B\), joining \(a\in A\) to \(b\in B\) whenever \(a\mid b\). Let \(F(A,x)\) denote the size of a maximum matching in \(G(A,x)\). Then
\[
f(m)=\min_{|A|=m}\ \min_{x\in\mathbb{R}} F(A,x),
\]
where the outer minimum is over all \(m\)-element sets \(A\subseteq \mathbb{N}\). This viewpoint will be especially relevant in the proof of the lower bound.

\subsection{AI-assisted proof discovery and formal verification}\label{sec:ai}

A distinctive feature of the present paper is the workflow by which the proof was obtained and verified. An initial proof draft was produced by ChatGPT (model: GPT-5.4 Pro),\footnote{GPT is a large language model developed by OpenAI; for background on the GPT-5 family, see \cite{OpenAIGPT5}.} which correctly identified the main proof strategy but left a gap in one of the detailed arguments. A public record of this stage of the process is available in the discussion thread for Erd\H{o}s Problem \#650; see \cite{Bloom650Thread}. The draft argument was then given to Aristotle,\footnote{Aristotle is an AI system for formal mathematical reasoning developed by Harmonic; see \cite{AristotlePaper}. For the project page, see \cite{AristotleProject}.} which supplied the missing details and produced a complete Lean formalization of the proof. For transparency, an early proof draft and a subsequent note identifying a gap in that draft are publicly available; see \cite{GPT650v1,GPT650Gap}. More information on the formalization process can be found in Section \ref{sec:lean}. 

This example suggests a possible research workflow in which a large language model is used primarily for proof search and high-level proof discovery, while a formal-reasoning system is used to turn the resulting argument into a fully rigorous and machine-checked proof. Such a workflow does not remove human judgment altogether (for example, the original proofs were overly lengthy and have, through human rewriting, decreased in size by about a factor of two), but it can substantially reduce the extent to which correctness depends on traditional line-by-line manual proof checking.

\subsection{Paper organization}

In Section~\ref{sec:statement} we state the main result and explain how the exact formula follows from separate upper and lower bounds. 
Section~\ref{sec:upper} proves the upper bound by constructing, for any product \(m=st\), a set $A$ of \(m\) integers and an interval of length \(2\max A\) in which at most \(s+t\) distinct multiples occur. 
Section~\ref{sec:lower} establishes the matching lower bound using a generalization of Hall's theorem, together with a careful analysis of the neighbourhood sizes in the relevant bipartite graph.  
Finally, Section~\ref{sec:lean} briefly records the process by which we obtained a formal verification of the argument in Lean.

\section{Statement and deduction of our main result}\label{sec:statement}

The main result we prove is the following.\footnote{Throughout, the symbol \doublecheck{} indicates that the proof of the corresponding statement has been formalized in Lean~4.}

\begin{vtheorem}\label{thm:main}
For every positive integer \(m\) we have
\[
f(m)=\min\bigl(m,\lceil 2\sqrt m\,\rceil\bigr).
\]
\end{vtheorem}

\begin{vremark} \label{rem:ge4}
As \(\lceil 2\sqrt{m}\,\rceil \le m\) for all \(m\ge 4\), Theorem~\ref{thm:main} in particular implies that \(f(m)=\lceil 2\sqrt{m}\,\rceil\) for all \(m\ge 4\).
\end{vremark}

Now, from the definition of \(f(m)\), the upper bound \(f(m)\le m\) is immediate. Accordingly, the proof of Theorem~\ref{thm:main} naturally splits into showing the upper bound $f(m) \le \lceil 2\sqrt m\,\rceil$ for all $m$ on the one hand, and the lower bound $f(m) \ge \min(m, \lceil 2\sqrt m\,\rceil)$ on the other. 

Given positive integers \(s\) and \(t\), for the upper bound we construct a specific set \(A\) with \(|A|=st\) and an open interval \(I\) of length \(2\max A\), such that \(I\) contains at most \(s+t\) distinct multiples of elements of \(A\). In particular, any matching in the corresponding bipartite graph has size at most \(s+t\), and hence
\begin{equation}\label{eq:fstlest_v1}
f(st)\le s+t.
\end{equation}
To deduce from this that \(f(m)\le \lceil 2\sqrt m\,\rceil\) for every \(m\), let \(k\) be such that \(k^2<m\le (k+1)^2\). If \(k^2<m\le\left(k+\tfrac12\right)^2\), then $m\le k(k+1)$ and \(\lceil 2\sqrt m\,\rceil=2k+1\). Since \(f(m)\) is monotone non-decreasing in \(m\) (indeed, one may always pass to a subset), we have \(f(m)\le f(k(k+1))\). Using the bound~\eqref{eq:fstlest_v1} with \(s=k\) and \(t=k+1\), we obtain
\[
f(m)\le f(k(k+1))\le 2k+1=\lceil 2\sqrt m\,\rceil.
\]
If instead \(\left(k+\tfrac12\right)^2 < m\le (k+1)^2\), then \(\lceil 2\sqrt m\,\rceil=2k+2\). In this case we take \(s=t=k+1\) in~\eqref{eq:fstlest_v1}, and similarly deduce
\[
f(m)\le f((k+1)^2)\le 2k+2=\lceil 2\sqrt m\,\rceil.
\]

We remark that \eqref{eq:fstlest_v1} is ever so slightly stronger than the earlier bound \(f(m^2)\le 2m\) proved by Erd\H{o}s and Selfridge~\cite{Erdos78,Erdos86}. Indeed, their bound yields \(f(m)\le 2\lceil \sqrt m\,\rceil\), which in general is weaker than \(f(m)\le \lceil 2\sqrt m\,\rceil\). Even though these estimates differ by at most $1$, this difference does of course become relevant if one wants to determine $f(m)$ precisely.

As for the lower bound, we create the bipartite graph that we introduced in Section \ref{sec:intro}. We then recall that, for a subset $S \subseteq A$, the neighbourhood $\Gamma(S) \subseteq B$ of $S$ is defined as the set of vertices in $B$ that are joined by an edge to an element in $S$. We then have the following generalization of Hall's theorem (see \cite[Exercise~16.2.8(b)]{BondyMurty}):

\begin{vlemma}\label{lem:deficientHall}
For any finite bipartite graph with vertex sets $A$ and $B$, there exists a matching of size $$\min \left(|A|, |A| - \max_{\emptyset \neq S \subseteq A}(|S| - |\Gamma(S)|) \right),$$ where the inner maximum is taken over all non-empty subsets $S \subseteq A$. 
\end{vlemma}

\begin{remark}
As in \cite{BondyMurty}, the above expression is also sometimes referred to as the K\"onig--Ore formula, and is generally more simply written as $$|A| - \max_{S \subseteq A}(|S| - |\Gamma(S)|)$$ instead, where the maximum is taken over \emph{all} subsets $S \subseteq A$. This is quickly seen to be equivalent to our version, by considering the empty set separately.
\end{remark}

With Lemma \ref{lem:deficientHall}, we claim that the desired inequality then follows from showing that for all $S \subseteq A$ we have
\begin{equation} \label{eq:gammabound}
|\Gamma(S)| \ge 2\sqrt{|S|}.
\end{equation}
To see why this would be sufficient, note that $t - 2\sqrt t$ is an increasing function on $\mathbb{N}$, due to the fact that the inequality $2 \sqrt{t+1} - 2\sqrt{t} < 1$ holds for all $t \ge 1$. With $|A| = m$, inequality \eqref{eq:gammabound} therefore implies $$|S| - |\Gamma(S)| \le |S| - 2\sqrt{|S|} \le m - 2\sqrt{m}$$ for all non-empty $S \subseteq A$, so that by Lemma \ref{lem:deficientHall} we deduce the existence of a matching of size at least $\min(m, 2 \sqrt{m})$. As $f(m)$ is an integer, this indeed implies $f(m) \ge \min(m, \left \lceil 2 \sqrt{m} \, \right \rceil)$. We will accomplish the proof of \eqref{eq:gammabound} by first splitting any interval of length $2 \max A$ into two halves. We then define an injection from $S$ into the product of the two separate halves of $\Gamma(S)$, in order to obtain the required lower bound on $|\Gamma(S)|$.

\section{The upper bound}\label{sec:upper}
As explained in Section \ref{sec:statement}, for the upper bound it is sufficient to show the following two-parameter inequality.

\begin{vtheorem}\label{thm:upper-st}
We have $f(st) \leq s+t$ for all positive integers $s$ and $t$.
\end{vtheorem}

\begin{proof}
If $\min(s, t) = 1$, then the upper bound $f(st) \leq st < s+t$ trivially holds. We will therefore assume $s, t \geq 2$, and by symmetry we may further assume $s \le t$. As mentioned in Section \ref{sec:statement}, we now aim to construct a set $A$ with $st$ elements for which an interval of length $2 \max A$ exists that contains at most $s+t$ multiples of the elements in $A$. 

In order to construct $A$, we let $D$ be the least common multiple of the integers $1, 2, \ldots, st$, and we further define
\[
\mathcal P:=
  \Bigl\{
  p>st:\ p \text{ is prime and } p\mid (q+rD)
  \text{ for some } 0 < |q| < s,\ |r| < t
  \Bigr\}.
\]
Since there are only finitely many pairs $(q,r)$ with $0 < |q| < s$ and $|r| < t$, and for each such pair the integer $q+rD$ is non-zero and hence has only finitely many prime divisors, the set $\mathcal P$ is finite.

Now, for each $p \in \mathcal P$, consider the residue classes $-i-jD \pmod p$, where $i$ and $j$ run over the integers with $1 \le i \le s$ and $1 \le j \le t$ respectively. In particular, the number of residue classes $-i-jD \pmod p$ is at most $st$. As $p > st$, these classes do not exhaust all residue classes modulo $p$, so that by the Chinese Remainder Theorem we deduce that there exists an integer $M > 2s + 2tD$ with $M \not \equiv -i-jD \pmod p$ for all $p \in \mathcal{P}$ and all $i$ and $j$. With this $M$, we define
\[
\alpha_{i,j}:= M+i+jD \qquad\text{and} \qquad A:= \{\alpha_{i,j} : 1 \le i \le s,\ 1 \le j \le t\}.
\]
To see that the $\alpha_{i,j} \in A$ are distinct (so that $|A| = st$), we note that $\alpha_{i,j} = \alpha_{k, l}$ implies that $D$ divides $k-i$. As $|k-i| < s < D$, we then get $i = k$, and consequently $j = l$. We furthermore have the following claim.

\begin{vclaim}\label{claim:gcd}
For any indices $(i,j)$ and $(k,l)$ with
$1\le i,k\le s$ and $1\le j,l\le t$, one has
\[
\gcd(\alpha_{i,j},\alpha_{k,l})\mid (i-k).
\]
\end{vclaim}

\begin{proof}
Without loss of generality we may assume $i \neq k$, as every integer divides $0$. Now, with $g:=\gcd(\alpha_{i,j},\alpha_{k,l})$ and with $v_p(\cdot)$ denoting the $p$-adic valuation, we aim to prove that $v_p(g) \le v_p(i-k)$ for all prime divisors $p$ of $g$.

Hence, let $p$ be an arbitrary prime divisor of $g$, and define $q := i-k \neq 0$ and $r := j-l$. Since $p \mid \alpha_{i,j}$ and \(\alpha_{i,j}-\alpha_{k,l}=q+rD\), we have
\[
M\equiv -i-jD \pmod p
\qquad\text{and}\qquad
p\mid(q+rD).
\]
Because $0 < |q| < s$ and $|r| < t$, if $p > st$, then the definition of $\mathcal P$
would imply that $p \in \mathcal P$, contradicting the choice of $M$. Hence $p \le st$, which we claim implies that $v_p(D) > v_p(q)$.

As $p \le st$ and \(D=\operatorname{lcm}(1,2,\dots,st)\), we have $v_p(D) > 0$, so that the inequality $v_p(D) > v_p(q)$ certainly holds if $v_p(q) = 0$. On the other hand, if $v_p(q) \ge 1$, then we note that $q^2 \mid D$, as $q^2 < s^2 \le st$. So in this case we get
\[
v_p(D) \ge v_p(q^2) = 2v_p(q) > v_p(q),
\]
where the final inequality uses $q \neq 0$. Thus indeed $v_p(D) > v_p(q)$ in both cases, from which we see \(v_p(q+rD) = v_p(q)\). Since $g \mid(q+rD)$, we now obtain
\[
v_p(g) \le v_p(q+rD) = v_p(q).
\]
As this holds for every prime $p$, we deduce that $g$ does indeed divide $q = i-k$.
\end{proof}

By the generalized Chinese Remainder Theorem (see \cite[Theorem~3.12]{JonesJonesENT}), an integer $x_0$ exists with $x_0 \equiv i \pmod{\alpha_{i,j}}$ for all $\alpha_{i,j} \in A$, as long as, for all $i, j, k, l$, the difference of the residue $i-k$ is divisible by $\gcd(\alpha_{i,j}, \alpha_{k,l})$. As this is precisely the content of Claim \ref{claim:gcd}, such an $x_0$ does indeed exist. 

We now define $x:=x_0-M$ and $I:=(x,x+2\max A)$, and claim that for every fixed pair $(i, j)$, there are at most two multiples of $\alpha_{i,j}$ in $I$. To see this, first note that the integer $x_0-i$ is a multiple of $\alpha_{i,j}$, by construction of $x_0$. As the multiples of $\alpha_{i,j}$ form an arithmetic progression with
common difference $\alpha_{i,j}$, it suffices to show that
\[
(x_0-i)-\alpha_{i,j}<x
\qquad\text{and}\qquad
(x_0-i)+2\alpha_{i,j}>x+2\max A.
\]
Indeed, on the one hand,
\begin{align*}
(x_0-i)-\alpha_{i,j}
&=x+M-i-(M+i+jD)\\
&=x-2i-jD\\
&<x.
\end{align*}
On the other hand,
\begin{align*}
(x_0-i)+2\alpha_{i,j}
&=x+M-i+2(M+i+jD)\\
&=x+3M+i+2jD\\
&>x+3M\\
&>x+2M+2s+2tD\\
&=x+2\max A.
\end{align*}
Any multiple of $\alpha_{i,j}$ lying in $I$ must therefore be either
\[
x_0-i
\qquad\text{or}\qquad
x_0-i+\alpha_{i,j}=x_0+M+jD.
\]
In particular, every multiple of an element of $A$ that lies in $I$ belongs to
\[
\{x_0-i:1\le i\le s\}\cup \{x_0+M+jD:1\le j\le t\}.
\]
As this latter set has at most $s+t$ elements, this finishes the proof.
\end{proof}

\begin{vremark} \label{rem:3max}
Let $\epsilon \in (0, 1)$ be arbitrary and recall that we chose $M > 2s + 2tD$. By instead choosing $M$ larger than $\epsilon^{-1} (3 - \epsilon) (s + tD)$, the above proof works even if we increase $I$ to $(x, x + (3 - \epsilon) \max A)$.
\end{vremark}

\section{The lower bound}\label{sec:lower}

We now prove the lower bound for $f(m)$.

\begin{vtheorem}\label{thm:lower}
For every positive integer $m$,
\[
f(m)\ge \min(m, \lceil 2\sqrt m\,\rceil).
\]
\end{vtheorem}

\begin{proof}
Fix \(m \in \mathbb{N}\), and let $A = \{a_1 < \cdots < a_m\}$ be any set of $m$ positive integers. Furthermore, with $x$ an arbitrary real number, define the interval $I := (x, x + 2a_m)$. As explained in Section \ref{sec:statement}, with $B:=I\cap\ZZ$ and $G$ the corresponding bipartite graph with vertex sets $A$ and $B$, it suffices to show inequality \eqref{eq:gammabound} for a fixed but arbitrary subset $S \subseteq A$. In order to do this, we need to distinguish between the case where $x$ is an integer multiple of $a_m$, and the case where it is not.

\medskip
\noindent
\underline{\emph{Case 1: $x \notin a_m\ZZ$.}}
We consider the partition \(B=B_-\sqcup B_+\) with
\[
B_-:=(x,x+a_m]\cap\ZZ,
\qquad \text{and} \qquad
B_+:=(x+a_m,x+2a_m)\cap\ZZ,
\]
and define
\[
\Gamma_-(S):=\Gamma(S)\cap B_-,
\qquad
\Gamma_+(S):=\Gamma(S)\cap B_+.
\]
Since $\Gamma_-(S)$ and $\Gamma_+(S)$ together partition $\Gamma(S)$, we note that
\begin{equation}\label{eq:gamma-sum}
|\Gamma(S)| = |\Gamma_-(S)| + |\Gamma_+(S)|.
\end{equation}
Moreover, we claim the inequality
\begin{equation}\label{eq:gamma-prod-case1}
|S|\le |\Gamma_-(S)|\,|\Gamma_+(S)|.
\end{equation}
To see this, for each $a \in S$, let $u_a$ be the largest multiple of $a$ that is contained in $B_-$. As $a$ is at most $a_m$, $u_a$ certainly exists. Moreover, by definition of $u_a$ and the assumption that $a_m$ does not divide $x$, we have that $u_a + a$ is a multiple of $a$ contained in $B_+$. Now consider the function
\[
\phi:S\to \Gamma_-(S)\times \Gamma_+(S),
\qquad
\phi(a):=(u_a,u_a+a).
\]
Since the difference of the two coordinates is exactly \(a\), this map is injective, which indeed gives \eqref{eq:gamma-prod-case1}. Combining \eqref{eq:gamma-sum} and \eqref{eq:gamma-prod-case1} with the AM–GM inequality already finishes the proof for this case, as we then get
\begin{align*}
|S| &\le |\Gamma_-(S)|\,|\Gamma_+(S)| \\
& \le \frac{\left(|\Gamma_-(S)| + |\Gamma_+(S)|\right)^2}{4} \\
&= \frac{|\Gamma(S)|^2}{4},
\end{align*}
implying \eqref{eq:gammabound}.

\medskip
\noindent
\underline{\emph{Case 2: $x \in a_m\ZZ$.}}
In this case we define $b_0 := x + a_m$, and let $G_0$ be the subgraph of $G$ with vertex sets $A_0 := A \setminus \{a_m \}$ and $B_0 := B \setminus \{b_0 \}$. Let $B_-$ and $B_+$ be defined as before, and define $\widetilde{B}_- := B_- \setminus \{b_0 \}$. For \(S\subseteq A_0\), we now write \(\Gamma_0(S)\) for its neighbourhood in \(G_0\), and we set
\[
\widetilde{\Gamma}_-(S):=\Gamma_0(S)\cap \widetilde{B}_-
\qquad \text{and} \qquad
\widetilde{\Gamma}_+(S):=\Gamma_0(S)\cap B_+.
\]
In analogy with \emph{Case 1}, it suffices to find an injective function $\phi:S \to \widetilde{\Gamma}_-(S) \times \widetilde{\Gamma}_+(S)$ for every subset $S \subseteq A_0$. Indeed, such a function would imply the analogous version of \eqref{eq:gamma-prod-case1}, from which $|\Gamma_0(S)| \ge 2\sqrt{|S|}$ follows. We would then by Lemma \ref{lem:deficientHall} deduce the existence of a matching of size at least $\min(m-1, 2\sqrt{m-1})$ in $G_0$, and adding the edge \(a_m \sim b_0\) finishes the proof, by providing a matching of size at least $$\min(m-1, 2\sqrt{m-1}) + 1 \ge \min(m, 2\sqrt m)$$ in the original graph $G$.

Hence, we fix a subset $S \subseteq A_0$ and, for each \(a\in S\), we let $u_a$ once again be the largest multiple of $a$ contained in $B_-$. If $a \nmid b_0$ we say that $a$ is of type (T1), and we define $\phi(a) = (u_a, u_a + a)$ as before. If $a$ does divide $b_0$, we remark that $u_a = b_0$. In this case, if we furthermore have $2a < a_m$ and $2a \nmid b_0$, then we say that $a$ is of type (T2), and we define $\phi(a) = (u_a - 2a, u_a + 2a)$. Finally, if $a \mid b_0$ but either $2a \ge a_m$ or $2a \mid b_0$, then we refer to $a$ as being of type (T3), and we define $\phi(a) = (u_a - a, u_a + a)$. With these definitions it is quickly verified that $\phi(S) \subseteq \widetilde{\Gamma}_-(S) \times \widetilde{\Gamma}_+(S)$, so it suffices to show that $\phi$ is injective.

For a start, as the difference of the coordinates is $a$, $4a$, or $2a$ (for elements $a$ of type (T1), (T2), (T3) respectively), we certainly have $\phi(a) \neq \phi(a')$ if distinct $a, a' \in S$ are of the same type. We will therefore by contradiction assume that $\phi(a) = \phi(a')$ with $a$ and $a'$ of different types, and note that this implies that at least one is a divisor of $b_0$, so that $b_0$ is exactly the average of the two coordinates.

First suppose that \(a\) is of type (T1) and \(a'\) is of type (T2). Then $\phi(a) = \phi(a')$ implies $a = 4a'$ by the difference of their coordinates. In particular, as $u_a$ is a multiple of $a$, $u_a$ must be divisible by $2a'$. However, this gives $$b_0 = u_a + \frac{a}{2} = 2a' \left(\frac{u_a}{2a'}\right) + 2a' = 2a'\left(\frac{u_a}{2a'} + 1\right),$$ contradicting the fact that $2a' \nmid b_0$ as \(a'\) is of type (T2).

Next, suppose that \(a\) is of type (T1) and \(a'\) is of type (T3), in which case $\phi(a) = \phi(a')$ implies $a = 2a'$. Since $2a' = a < a_m$, this implies $2a' \mid b_0$ by the assumption that $a'$ is of type (T3). This also leads to a contradiction, as $$b_0 = u_a + \frac{a}{2} = 2a' \cdot \frac{u_a}{2a'} + a' = a'\left(2 \cdot \frac{u_a}{2a'} + 1\right),$$ which is an odd multiple of $a'$.

Finally, suppose that \(a\) is of type (T2) and \(a'\) is of type (T3). Then $\phi(a) = \phi(a')$ implies $a' = 2a$, while $a' \mid b_0$ as $a'$ is of type (T3). But this contradicts the fact that $2a \nmid b_0$. 

The opposite orderings are symmetric, so no mixed-type collision can occur. Hence $\phi$ is indeed injective, finishing the proof.
\end{proof}

\section{Lean formalization}\label{sec:lean}

We used Aristotle, the automated theorem proving tool from Harmonic, in order to obtain formal proofs of all results in this paper. One interesting feature of this formalization process was that we discovered, only after the formalization had already been completed, that the original proof draft produced by ChatGPT contained a gap. That gap was subsequently repaired, and the present paper reflects the corrected argument. More precisely, on closer inspection it became clear that ChatGPT's initial proposal for $\phi:S \to \widetilde{\Gamma}_-(S) \times \widetilde{\Gamma}_+(S)$ defined in \emph{Case 2} of the lower bound, was not an injection. This was of course worrying, but also confusing, as the formal proof had already been finished.

We then had to translate the Lean proof back into informal language, to figure out what had happened. It then turned out that Aristotle had not only noticed the gap, but also managed to come up with a variation of $\phi$ that actually did work! Hence, Aristotle's definition of $\phi$ is the one that we actually used in Section \ref{sec:lower} in order to prove Theorem \ref{thm:lower}.

As for the accompanying Lean file~\cite{LeanFile650}, the file header records the following metadata:
\begin{itemize}
  \item Lean version: \texttt{leanprover/lean4:v4.28.0};
  \item Mathlib version: \texttt{8f9d9cff6bd728b17a24e163c9402775d9e6a365};
  \item The formalization was obtained by Aristotle from Harmonic.
\end{itemize}

The file consists of $1116$ lines and contains $8$ definitions, $42$ lemmas, and $5$ named theorems:

\begin{itemize}
  \item \texttt{erdos\_f\_upper\_bound},
  \item \texttt{large\_3maxA\_version},
  \item \texttt{erdos\_f\_lower\_bound},
  \item \texttt{erdos\_f\_eq}.
	\item \texttt{erdos\_f\_eq\_ge4}.
\end{itemize}

These latter $5$ results refer to Theorem \ref{thm:upper-st}, Remark \ref{rem:3max}, Theorem \ref{thm:lower}, Theorem \ref{thm:main} and Remark \ref{rem:ge4} respectively.

For background on Lean and Mathlib, see~\cite{Mathlib,Lean4}.

\end{document}